\documentclass{article}
\usepackage{amsmath,amsfonts,multirow}

\newcommand{\textlineskip}{\baselineskip=11.5dd plus0.5dd minus0.2dd}

\def\qed{\hbox{${\vcenter{\vbox{
   \hrule height 0.4pt\hbox{\vrule width 0.4pt height 6pt
   \kern5pt\vrule width 0.4pt}\hrule height 0.4pt}}}$}}

\def\tcaption#1{\addtocounter{table}{1} \begin{center} {\footnotesize
TABLE~\thetable. #1} \end{center}}

\catcode`\@=11

\font\twrm=cmr12

\font\eightrm=cmr8

\newcounter{subsect}

\def\section#1{\par \bigskip \setcounter{subsect}{0}
\setcounter{equation}{0}
\addtocounter{section}{1}\begin{center}\thesection.
\uppercase{#1} \end{center} \par \smallskip }

\def\subsection#1{\par \medskip
\addtocounter{subsect}{1}\begin{center}{ \eightrm
\thesection.\thesubsect. \uppercase{#1}} \end{center} \par
\smallskip }

\makeatletter
\newcommand*{\rom}[1]{\expandafter\@slowromancap\romannumeral #1@}
\makeatother
\newcommand{\Bia}{\operatorname{Bia}}

\newcommand{\ie}{i.e. }

\newcommand{\F}{\mathcal{F}}

\newcommand{\M}{(M,\f,\xi,\eta,g)}
\newcommand{\Lf}{(L,\f,\xi,\eta,g)}

\newcommand{\R}{\mathbb{R}}

\newcommand{\n}{\nabla}
\newcommand{\nn}{\tilde{\n}}
\newcommand{\f}{\varphi}
\newcommand{\tg}{\tilde{g}}
\newcommand{\om}{\omega}

\newcommand{\al}{\alpha}

\newcommand{\ta}{\theta}
\newcommand{\lm}{\lambda}

\newcommand{\norm}[1]{\left\Vert#1\right\Vert ^2}

\newcommand{\Id}{\mathrm{Id}}

\newcommand{\ff}{\varphi}

\begin{document}

\normalsize\textlineskip
\thispagestyle{empty}

 \begin{center} {\twrm

\vspace*{1.0in}
ALMOST CONTACT B-METRIC STRUCTURES\\
\vspace*{0.035truein}
AND THE BIANCHI CLASSIFICATION \\
\vspace*{0.035truein}
OF THE THREE-DIMENSIONAL LIE ALGEBRAS
}\\
\vspace*{2cc}
{\eightrm HRISTO M.~MANEV}
\end{center}

\vspace*{30dd}

{\parindent0pt \footnotesize \leftskip20pt
\rightskip20pt  \baselineskip10pt
The object of investigation are the almost contact manifolds with
B-metric in the lowest dimension three, constructed on Lie algebras.
It is considered a relation between the classes in the
Bianchi classification of three-dimensional real Lie algebras and the classes of
a classification of the considered manifolds. There are studied some
geometrical characteristics in some special classes.
\smallskip

{\bf  2010
Math.\ Subject Classification}.  53C15, 53C50, 53D15 
\par
}

\vspace*{16dd}

\section{Introduction}

The differential geometry of the manifolds equipped with an almost contact structure
is well studied (e.g. [3]
). The almost contact manifolds with B-metric are introduced
and classified in [6]
. These manifolds are the odd-dimensional
counterpart of the almost complex manifolds with Norden metric
[5, 7].

An object of special interest is the case of the lowest dimension
of the considered manifolds. We investigate the almost contact B-metric
manifolds in dimension three and get explicit results.
Some curvature identities of the
three-dimensional manifolds of this type are studied in
[11, 12]. 

Almost contact manifolds with B-metric can be constructed on Lie algebras.
It is known that all three-dimensional real Lie algebras
are classified in [1, 2]
. The main goal of this paper is to find
a relation between the classes in the Bianchi classification
and the classification of almost contact B-metric manifolds given in [6]. 
Moreover, the present work gives some geometrical characteristics of considered manifolds in certain special classes.

The present paper is organized as follows. In
Sect.~2 
we recall
some preliminary facts about the almost contact B-metric manifolds. In Sect.~3 
we equip each Bianchi-type Lie algebra with an almost contact B-metric structure. In Sect.~4 
we give the relation between the Bianchi classification and the classification given in [6]
. Sect.~5 
is devoted to the curvature properties of some of the considered manifolds.

\section{Preliminaries}\label{sect-prel}

Let $(M,\f,\xi,\eta,g)$ be an almost contact manifold with
B-met\-ric or an \emph{almost contact B-metric manifold}, where $M$
is a $(2n+1)$-dimensional differentiable manifold, $(\f,\xi,\eta)$ is an almost
contact structure consisting of an endomorphism
$\f$ of the tangent bundle, a Reeb vector field $\xi$ and its dual contact 1-form
$\eta$. Moreover, $M$ is equipped with a pseu\-do-Rie\-mannian
metric $g$, called a \emph{B-metric}, such that the following
algebraic relations are satisfied [6]:
\[
\begin{array}{c}
\f\xi = 0,\qquad \f^2 = -\Id + \eta \otimes \xi,\qquad
\eta\circ\f=0,\qquad \eta(\xi)=1,\\[0pt]
g(\f x, \f y) = - g(x,y) + \eta(x)\eta(y),
\end{array}
\]
where $\Id$ is the identity. In the latter equalities and further, $x$, $y$, $z$, $w$ will stand for arbitrary elements of the
algebra of the smooth vector fields on $M$ or vectors in the tangent space $T_pM$ of $M$ at an arbitrary
point $p$ in $M$.

The associated B-metric $\tg$ of $g$ is determined by
\(\tg(x,y)=g(x,\f y)+\eta(x)\eta(y)\).  The manifold
$(M,\f,\xi,\eta,\tg)$ is also an almost contact B-metric manifold.
The signature of both metrics $g$ and $\tg$ is necessarily
$(n+1,n)$. We denote the Levi-Civita connection
of $g$ and $\tg$ by $\n$ and $\nn$, respectively.

A classification of almost contact B-metric manifolds, consisting of eleven basic
classes $\F_1$, $\F_2$, $\dots$, $\F_{11}$, is given in
[6]. 
This classification is made with respect
to the tensor $F$ of type (0,3) defined by
\begin{equation}\label{F=nfi}
F(x,y,z)=g\bigl( \left( \nabla_x \f \right)y,z\bigr)
\end{equation}
and having the following properties:
\[
F(x,y,z)=F(x,z,y)=F(x,\f y,\f z)+\eta(y)F(x,\xi,z)
+\eta(z)F(x,y,\xi).
\]

The special class determined by the condition $F(x,y,z)=0$ is denoted by $\F_0$.
This class is the intersection of all the basic classes.
Hence $\F_0$ is the
class of almost contact B-metric manifolds with $\n$-parallel
structures, \ie $\n\f=\n\xi=\n\eta=\n g=\n \tg=0$. Therefore $\F_0$ is the class of the \emph{cosymplectic manifolds with B-metric}.

According to [10], 
the \emph{square norm of $\nabla \f$} is defined by:
\begin{equation}\label{mi:snf}
    \norm{\nabla \ff}=g^{ij}g^{ks}
    g\bigl(\left(\nabla_{e_i} \ff\right)e_k,\left(\nabla_{e_j}
    \ff\right)e_s\bigr).
\end{equation}
It is clear, $\norm{\nabla \ff}=0$ is valid if $\M$
is a cosymplectic manifold with B-metric, but the inverse implication is not always
true. An almost contact B-metric manifold having a zero square
norm of $\n\ff$ is called an
\emph{isotropic-cosymplectic B-metric manifold}.

If $\left\{e_i;\xi\right\}$ $(i=1,2,\dots,2n)$ is a basis of
$T_pM$ and $\left(g^{ij}\right)$ is the inverse matrix of
$\left(g_{ij}\right)$, then the 1-forms
$\theta$, $\theta^*$, $\omega$, called \emph{Lee forms}, are associated with $F$ and defined
by:
\[
\theta(z)=g^{ij}F(e_i,e_j,z),\quad \theta^*(z)=g^{ij}F(e_i,\f
e_j,z), \quad \omega(z)=F(\xi,\xi,z).
\]

Let now consider the case of the lowest dimension of the almost contact B-metric manifold
$M$, \ie $\dim{M}=3$.

We introduce an almost contact structure $(\f,\xi,\eta)$ on $M$ defined by
\begin{equation}\label{strL}
\begin{array}{c}
\f e_1=e_{2},\quad \f e_{2}=- e_1,\quad \f e_3=0,\quad \xi=
e_3,\quad \\[0pt]
\eta(e_1)=\eta(e_{2})=0,\quad \eta(e_3)=1
\end{array}
\end{equation}
and a B-metric $g$ such that
\begin{equation}\label{gij}
g(e_1,e_1)=-g(e_2,e_2)=g(e_3,e_3)=1,\quad g(e_i,e_j)=0,\;\; i\neq j
\in \{1,2,3\}.
\end{equation}

Let us denote the components $F_{ijk}=F(e_i,e_j,e_k)$ of $F$ with respect to a \emph{$\f$-basis}
$\left\{e_1,e_2,e_3\right\}$ of $T_pM$.

According to [8], 
the components of the Lee forms are
\[
\begin{array}{lll}
\ta_1=F_{111}-F_{221},\qquad &\ta_2=F_{112}-F_{211},\qquad &\ta_3=F_{113}-F_{223},\\[0pt]
\ta^*_1=F_{112}+F_{211},\qquad &\ta^*_2=F_{111}+F_{221},\qquad &\ta^*_3=F_{123}+F_{213},\\[0pt]
\om_1=F_{331},\qquad &\om_2=F_{332},\qquad &\om_3=0.
\end{array}
\]

Then, if $F_s$ $(s=1,2,\dots,11)$ are the components of $F$ in the
corresponding basic classes $\F_s$ and $x=x^ie_i$, $y=y^je_j$,
$z=z^ke_k$ for arbitrary vectors in $T_pM$, we have [8]:
\begin{equation}\label{Fi3}
\begin{array}{l}
F_{1}(x,y,z)=\left(x^1\ta_1-x^2\ta_2\right)\left(y^1z^1+y^2z^2\right),\\[0pt]
\qquad\ta_1=F_{111}=F_{122},\qquad \ta_2=-F_{211}=-F_{222}; \\[0pt]
F_{2}(x,y,z)=F_{3}(x,y,z)=0;
\\[0pt]
F_{4}(x,y,z)=\frac{1}{2}\ta_3\Bigl\{x^1\left(y^3z^1+y^1z^3\right)
-x^2\left(y^3z^2+y^2z^3\right)\bigr\},\\[0pt]
\qquad \frac{1}{2}\ta_3=F_{131}=F_{113}=-F_{232}=-F_{223};\\[0pt]
F_{5}(x,y,z)=\frac{1}{2}\ta^*_3\bigl\{x^1\left(y^3z^2+y^2z^3\right)
+x^2\left(y^3z^1+y^1z^3\right)\bigr\},\\[0pt]
\qquad \frac{1}{2}\ta^*_3=F_{132}=F_{123}=F_{231}=F_{213};\\[0pt]
F_{6}(x,y,z)=F_{7}(x,y,z)=0;\\[0pt]
F_{8}(x,y,z)=\lm\bigl\{x^1\left(y^3z^1+y^1z^3\right)
+x^2\left(y^3z^2+y^2z^3\right)\bigr\},\\[0pt]
\qquad \lm=F_{131}=F_{113}=F_{232}=F_{223};\\[0pt]
F_{9}(x,y,z)=\mu\bigl\{x^1\left(y^3z^2+y^2z^3\right)
-x^2\left(y^3z^1+y^1z^3\right)\bigr\},\\[0pt]
\qquad \mu=F_{132}=F_{123}=-F_{231}=-F_{213};\\[0pt]
F_{10}(x,y,z)=\nu x^3\left(y^1z^1+y^2z^2\right),\qquad
\nu=F_{311}=F_{322};\\[0pt]
F_{11}(x,y,z)=x^3\bigl\{\left(y^1z^3+y^3z^1\right)\om_{1}
+\left(y^2z^3+y^3z^2\right)\om_{2}\bigr\},\\[0pt]
\qquad \om_1=F_{313}=F_{331},\qquad \om_2=F_{323}=F_{332}.
\end{array}
\end{equation}

Obviously, the class of three-dimensional almost contact B-metric
manifolds is
\[
\F_1 \oplus \F_4 \oplus \F_5 \oplus \F_8 \oplus \F_9 \oplus
\F_{10} \oplus \F_{11}.
\]

Let $R=\left[\n,\n\right]-\n_{[\ ,\ ]}$ be the
curvature (1,3)-tensor of $\nabla$. The corresponding curvature
$(0,4)$-tensor is denoted by the same letter: $R(x,y,z,w)$
$=g(R(x,y)z,w)$. The following properties are valid:
\[
\begin{array}{c}
    R(x,y,z,w)=-R(y,x,z,w)=-R(x,y,w,z), \\[0pt]
R(x,y,z,w)+R(y,z,x,w)+R(z,x,y,w)=0.
\end{array}
\]

It is known from [11] 
that every 3-dimen\-sion\-al cosymplectic B-metric manifold is flat, \ie $R
= 0$.

 The Ricci
tensor $\rho$ and the scalar curvature $\tau$ for $R$ as well as
their associated quantities are defined respectively by
\[
\begin{array}{ll}
    \rho(y,z)=g^{ij}R(e_i,y,z,e_j),\qquad &
    \tau=g^{ij}\rho(e_i,e_j),\\[0pt]
    \rho^*(y,z)=g^{ij}R(e_i,y,z,\f e_j),\qquad &
    \tau^*=g^{ij}\rho^*(e_i,e_j),
\end{array}
\]
where $\left\{e_1,e_2,\dots,e_{2n+1}\right\}$ is an arbitrary basis of $T_pM$.

Let $\al$ be a non-degenerate 2-plane (section) in $T_pM$.
It is known that the special 2-planes with respect to $(\f,\xi,\eta,g)$ are:
a \emph{totally real section} if $\al$ is orthogonal to its
$\f$-image $\f\al$, a
\emph{$\f$-holomorphic  section} if $\al$ coincides with
$\f\al$ and a
\emph{$\xi$-section} if $\xi$ lies on $\al$.

The sectional curvature $k(\al; p)(R)$ of
$\al$ with an arbitrary basis $\{x,y\}$ at $p$
is defined by
\[
k(\al; p)(R)=\frac{R(x,y,y,x)}{g(x,x)g(y,y)-g(x,y)^2}.%
\]

According to [9], 
it is reasonable to
call a manifold $M$ whose Ricci tensor satisfies the condition
\[
\rho=\lm g + \mu \tg + \nu \eta\otimes\eta
\]
 an
\emph{$\eta$-complex-Einstein manifold}.

\section{Equipping of each Bianchi-type Lie algebra with almost contact B-metric structure}\label{sect-2}
It is known that L. Bianchi has categorized all
three-dimensional real (and complex) Lie algebras. He proved that every
three-dimensional Lie algebra is isomorphic to one, and only one, Lie algebra
of his list (cf. [1, 2]
). These isomorphism classes form the so-called Bianchi classification and
are noted by $\Bia$(I), $\Bia$(II), $\Bia$(IV), $\Bia$(V), $\Bia$(VI$_h$)\ ($h \leq 0$), $\Bia$(VII$_h$)\
($h \geq 0$), $\Bia$(VIII) and $\Bia$(IX). The class $\Bia$(III) coincides with $\Bia$(VI$_{-1}$). The following theorem introduces the Bianchi classification.

{\bf Theorem~A}.
([1, 2]) 
Let $\mathfrak{l}$ be a real three-dimensional Lie algebra.
Then $\mathfrak{l}$ is isomorphic to exactly one of the following Lie algebras $(\R^3, [\cdot,\cdot])$,
where the Lie bracket is given on the canonical basis $\{e_1, e_2, e_3\}$ as follows:
\[
\begin{array}{llll}
\Bia(\mathrm{I}):& [e_1, e_2] = o,& [e_2, e_3] = o,& [e_3, e_1] = o;\\
\Bia(\mathrm{II}):& [e_1, e_2] = o,& [e_2, e_3] = e_1,& [e_3, e_1] = o;\\
\Bia(\mathrm{IV}) :& [e_1, e_2] = o,& [e_2, e_3] = e_1 - e_2,& [e_3, e_1] = e_1;\\
\Bia(\mathrm{V}) :& [e_1, e_2] = o,& [e_2, e_3] = e_2,& [e_3, e_1] = e_1;\\
\Bia(\mathrm{VI}_h)\ (h \leq 0):& [e_1, e_2] = o,& [e_2, e_3] = e_1 - he_2,& [e_3, e_1] = he_1 - e_2;\\
\Bia(\mathrm{VII}_h)\ (h \geq 0):& [e_1, e_2] = o,& [e_2, e_3] = e_1 - he_2,& [e_3, e_1] = he_1 + e_2;\\
\Bia(\mathrm{VIII}):& [e_1, e_2] = -e_3,& [e_2, e_3] = e_1,& [e_3, e_1] = e_2;\\
\Bia(\mathrm{IX}):& [e_1, e_2] = e_3,& [e_2, e_3] = e_1,& [e_3, e_1] = e_2,\\
\end{array}
\]
where $o$ is the zero vector of $\mathfrak{l}$.
\vspace*{12pt}

The geometrization conjecture, associated with W. Thurston, states that every closed manifold of dimension three could be decomposed in a canonical way into pieces, connected to one of eight types of Thurston's geometric structures ([13]): 
Euclidean geometry $E^3$, Spherical geometry $S^3$, Hyperbolic geometry $H^3$, the geometry of $S^2\times\R$, the geometry of $H^2\times\R$, the geometry of the universal cover $\widetilde{SL}(2,\R)$ of the special linear group $SL(2,\R)$, the $Nil$ geometry, the $Solv$ geometry.

Seven of the eight Thurston geometries can be associated to a class of the Bianchi classification as it is shown in the following table. 
The Thurston geometry on $S^2\times\R$ has no such a realization (e.g. [4]).

{\footnotesize{
\tcaption{Relations between the Bianchi types and the Thurston geometries}  \vskip 5pt
\begin{center}
\begin{tabular}{|l|l|}
  \hline
  $\Bia(\mathrm{I})$ & $E^3$    \\
  \hline
  $\Bia(\mathrm{II})$ & $Nil$  \\
    \hline
  $\Bia(\mathrm{III})$ & $H^2\times\R$ \\
    \hline
    $\Bia(\mathrm{IV})$ &  \\
    \hline
    $\Bia(\mathrm{V})$ & $H^3$  \\
    \hline
    $\Bia(\mathrm{VI}_0)$ & $Solv$ \\
    \hline
\end{tabular}
\quad
\begin{tabular}{|l|l|}
  \hline
    $\Bia(\mathrm{VI}_{h<0})$ & \\
    \hline
    $\Bia(\mathrm{VII}_0)$ & $E^3$  \\
    \hline
    $\Bia(\mathrm{VII}_{h>0})$ & \\
    \hline
    $\Bia(\mathrm{VIII})$ & $\widetilde{SL}(2,\R)$  \\
    \hline
    $\Bia(\mathrm{IX})$ & $S^3$ \\
  \hline
\end{tabular}
\end{center}
}}

\vspace*{12pt}
Let us consider each Lie algebra from the Bianchi classification,
equipped with an almost contact structure $(\f,\xi,\eta)$ and a B-metric $g$ as in \eqref{strL} and \eqref{gij}.

The presence of the structure $(\f,\xi,\eta, g)$ gives us a reason to consider the relation between the Bianchi types and the classification of almost contact B-metric manifolds in [6]. 

We obtain immediately the following
\vspace*{12pt}

{\bf Proposition~3.1}.
Some Bianchi types can be equipped with a structure $(\f,\xi,\eta,\allowbreak{} g)$ in several ways. 
In the cases $\Bia(\mathrm{I})$ and $\Bia(\mathrm{IX})$ there are only one variant.
In the rest cases, there are three possible subtypes of each type
differing each other by a cyclic change of the basic vectors $e_1$, $e_2$, $e_3$.
All subtypes are determined as follows:

{\footnotesize{
\tcaption{Equipping of the Bianchi types Lie algebras with a $(\f,\xi,\eta,\allowbreak{} g)$ structure}  \vskip 5pt
\begin{center}
\begin{tabular}{|l|lll|}
  \hline
  \multicolumn{4}{|l|}{$\Bia(\mathrm{I})$}\\\hline
  (1) & $[e_1,e_2]=o$, \, & $[e_2,e_3]=o$,\, & $[e_3,e_1]=o$ \\
  \hline
  \multicolumn{4}{|l|}{$\Bia(\mathrm{II})$}\\\hline
  (1) & $[e_1,e_2]=o$,\, &   $[e_2,e_3]=e_1$,\, & $[e_3,e_1]=o$ \\
  (2) & $[e_1,e_2]=o$,\, &   $[e_2,e_3]=o$,\, &   $[e_3,e_1]=e_2$ \\
  (3) & $[e_1,e_2]=e_3$,\, & $[e_2,e_3]=o$,\, &   $[e_3,e_1]=o$ \\
  \hline
  \multicolumn{4}{|l|}{$\Bia(\mathrm{III})\equiv \Bia(\mathrm{VI}_{-1})$}\\\hline
  (1) & $[e_1,e_2]=o$,\, &           $[e_2,e_3]=e_1+e_2$,\, &     $[e_3,e_1]=-e_1-e_2$ \\
  (2) & $[e_1,e_2]=-e_2-e_3$,\, &    $[e_2,e_3]=o$,\, &           $[e_3,e_1]=e_2+e_3$ \\
  (3) & $[e_1,e_2]=e_1+e_3$,\, &     $[e_2,e_3]=-e_1-e_3$,\, &    $[e_3,e_1]=o$ \\
  \hline
  \multicolumn{4}{|l|}{$\Bia(\mathrm{IV})$}\\\hline
  (1) & $[e_1,e_2]=o$,\, &           $[e_2,e_3]=e_1-e_2$,\, &         $[e_3,e_1]=e_1$ \\
  (2) & $[e_1,e_2]=e_2$,\, &         $[e_2,e_3]=o$,\, &               $[e_3,e_1]=e_2-e_3$ \\
  (3) & $[e_1,e_2]=-e_1+e_3$,\, &    $[e_2,e_3]=e_3$,\, &             $[e_3,e_1]=o$ \\
  \hline
  \multicolumn{4}{|l|}{$\Bia(\mathrm{V})$}\\\hline
  (1) & $[e_1,e_2]=o$,\, &       $[e_2,e_3]=e_2$,\, &     $[e_3,e_1]=e_1$ \\
  (2) & $[e_1,e_2]=e_2$,\, &     $[e_2,e_3]=o$,\, &       $[e_3,e_1]=e_3$ \\
  (3) & $[e_1,e_2]=e_1$,\, &     $[e_2,e_3]=e_3$,\, &     $[e_3,e_1]=o$ \\
  \hline
  \multicolumn{4}{|l|}{$\Bia(\mathrm{VI}_{h})$, $h\leq0$}\\\hline
  (1) & $[e_1,e_2]=o$,\, &           $[e_2,e_3]=e_1-he_2$,\, &    $[e_3,e_1]=he_1-e_2$ \\
  (2) & $[e_1,e_2]=he_2-e_3$,\, &    $[e_2,e_3]=o$,\, &           $[e_3,e_1]=e_2-he_3$ \\
  (3) & $[e_1,e_2]=-he_1+e_3$,\, &   $[e_2,e_3]=-e_1+he_3$,\, &   $[e_3,e_1]=o$ \\
  \hline
  \multicolumn{4}{|l|}{$\Bia(\mathrm{VII}_{h})$, $h\geq0$}\\\hline
  (1) & $[e_1,e_2]=o$,\, &           $[e_2,e_3]=e_1-he_2$,\, &    $[e_3,e_1]=he_1+e_2$ \\
  (2) & $[e_1,e_2]=he_2+e_3$,\, &    $[e_2,e_3]=o$,\, &           $[e_3,e_1]=e_2-he_3$ \\
  (3) & $[e_1,e_2]=-he_1+e_3$,\, &   $[e_2,e_3]=e_1+he_3$,\, &    $[e_3,e_1]=o$ \\
  \hline
  \multicolumn{4}{|l|}{$\Bia(\mathrm{VIII})$}\\\hline
  (1) & $[e_1,e_2]=-e_3$,\, &    $[e_2,e_3]=e_1$,\, &     $[e_3,e_1]=e_2$ \\
  (2) & $[e_1,e_2]=e_3$,\, &     $[e_2,e_3]=-e_1$,\, &    $[e_3,e_1]=e_2$ \\
  (3) & $[e_1,e_2]=e_3$,\, &     $[e_2,e_3]=e_1$,\, &     $[e_3,e_1]=-e_2$ \\
  \hline
  \multicolumn{4}{|l|}{$\Bia(\mathrm{IX})$}\\\hline
  (1) & $[e_1,e_2]=e_3$,\, & $[e_2,e_3]=e_1$,\, & $[e_3,e_1]=e_2$ \\
  \hline
\end{tabular}
\end{center}
}}

\section{Almost contact B-metric manifolds of each Bianchi type}\label{sect-3}
Let us consider the Lie group $L$ corresponding to the given Lie algebra $\mathfrak{l}$.
Each definition of a Lie algebra for the different subtypes in Proposition~3.1 
generates a corresponding almost contact B-metric manifold denoted by $\Lf$. In this section we characterize the obtained manifolds with respect to the classification in [6]. 

Using \eqref{Fi3}, we obtain the corresponding components of $F$ in each subtypes (1), (2), (3) in Proposition~3.1 
and determine the corresponding class of almost contact B-metric manifolds. The results are given in the following
\vspace*{12pt}

{\bf Theorem~4.1}. 
The manifold $\Lf$, determined by each type of Lie algebra given in Proposition~3.1, 
belongs to a class of the classification in $[6]$ 
 as it is given in the following table:

{\footnotesize{
\tcaption{Relations between the Bianchi types and the classes in $[6]$}  \vskip 5pt
\begin{center}
\begin{tabular}{|l|l|}
  \hline
  \multicolumn{2}{|l|}{$\Bia(\mathrm{I})$}\\\hline
  (1) & $\F_0$ \\
  \hline
  \multicolumn{2}{|l|}{$\Bia(\mathrm{II})$}\\\hline
  (1) & $\F_4\oplus\F_{10}$ \\
  (2) & $\F_4\oplus\F_{10}$ \\
  (3) & $\F_8\oplus\F_{10}$ \\
  \hline
  \multicolumn{2}{|l|}{$\Bia(\mathrm{III})$}\\\hline
  (1) & $\F_5\oplus\F_{10}$ \\
  (2) & $\F_1\oplus\F_4\oplus\F_8\oplus\F_{11}$ \\
  (3) & $\F_1\oplus\F_4\oplus\F_8\oplus\F_{10}\oplus\F_{11}$ \\
  \hline
  \multicolumn{2}{|l|}{$\Bia(\mathrm{IV})$}\\\hline
  (1) & $\F_4\oplus\F_5\oplus\F_{10}$ \\
  (2) & $\F_1\oplus\F_4\oplus\F_{10}\oplus\F_{11}$ \\
  (3) & $\F_1\oplus\F_8\oplus\F_{10}\oplus\F_{11}$ \\
  \hline
  \multicolumn{2}{|l|}{$\Bia(\mathrm{V})$}\\\hline
  (1) & $\F_9$ \\
  (2) & $\F_1\oplus\F_{11}$ \\
  (3) & $\F_1\oplus\F_{11}$ \\
  \hline
  \multicolumn{2}{|l|}{$\Bia(\mathrm{VI}_0)$}\\\hline
  (1) & $\F_{10}$ \\
  (2) & $\F_4\oplus\F_8$ \\
  (3) & $\F_4\oplus\F_8\oplus\F_{10}$ \\
  \hline
\end{tabular}\qquad
\begin{tabular}{|l|l|}
  \hline
  \multicolumn{2}{|l|}{$\Bia(\mathrm{VI}_{h})$, $h<0$}\\\hline
  (1) & $\F_5\oplus\F_{10}$ \\
  (2) & $\F_1\oplus\F_4\oplus\F_8\oplus\F_{11}$ \\
  (3) & $\F_1\oplus\F_4\oplus\F_8\oplus\F_{10}\oplus\F_{11}$ \\
  \hline
  \multicolumn{2}{|l|}{$\Bia(\mathrm{VII}_0)$}\\\hline
  (1) & $\F_4$ \\
  (2) & $\F_4\oplus\F_8\oplus\F_{10}$ \\
  (3) & $\F_4\oplus\F_8$ \\
  \hline
  \multicolumn{2}{|l|}{$\Bia(\mathrm{VII}_{h})$, $h>0$}\\\hline
  (1) & $\F_4\oplus\F_5$ \\
  (2) & $\F_1\oplus\F_4\oplus\F_8\oplus\F_{10}\oplus\F_{11}$ \\
  (3) & $\F_1\oplus\F_4\oplus\F_8\oplus\F_{11}$ \\
  \hline
  \multicolumn{2}{|l|}{$\Bia(\mathrm{VIII})$}\\\hline
  (1) & $\F_4\oplus\F_8\oplus\F_{10}$ \\
  (2) & $\F_8\oplus\F_{10}$ \\
  (3) & $\F_8\oplus\F_{10}$ \\
  \hline
  \multicolumn{2}{|l|}{$\Bia(\mathrm{IX})$}\\\hline
  (1) & $\F_4\oplus\F_8\oplus\F_{10}$ \\
  \hline
\end{tabular}
\end{center}
}}

\vspace*{12pt}
{\it Proof.}
We give our arguments for the case of $\Bia$(II). In a similar way we prove the other cases.

Using Theorem A, Eq. \eqref{gij} and the
Koszul equality
\[
2g\left(\n_{e_i}e_j,e_k\right)
=g\left([e_i,e_j],e_k\right)+g\left([e_k,e_i],e_j\right)
+g\left([e_k,e_j],e_i\right),
\]
we obtain the components of the Levi-Civita connection $\n$ of $g$.
Then, by them, \eqref{F=nfi} and \eqref{strL},
we get the following nonzero components $F_{ijk}$ and $\ta_k$ for the different subtypes:
\[
\begin{array}{l}
(1)\quad F_{113}=F_{131}=-F_{223}=-F_{232}=-\frac{1}{2}, \quad
   F_{311}=F_{322}=-1,\quad \ta_3=-1; \\
(2)\quad F_{113}=F_{131}=-F_{223}=-F_{232}=-\frac{1}{2}, \quad
 F_{311}=F_{322}=1,\quad  \ta_3=-1; \\
(3)\quad F_{113}=F_{131}=F_{223}=F_{232}=\frac{1}{2}, \quad
                  F_{311}=F_{322}=1. \\
\end{array}
\]
After that, bearing in mind \eqref{Fi3}, we conclude the corresponding class of each subtype of $\Bia$(II) as follows:
\[
\begin{array}{l}
(1)\quad \Lf \in  \F_4\oplus\F_{10};  \\
(2)\quad \Lf \in  \F_4\oplus\F_{10}; \\
(3)\quad \Lf  \in  \F_8\oplus\F_{10}. \ \qed
\end{array}
\]

\newpage

\section{Curvature properties of the considered manifolds in some Bianchi classes}\label{sect-4}

Now, we focuss our considerations on the Bianchi classes depending on real parameter $h$.
They are $\Bia$(VI$_h$) and $\Bia$(VII$_h$). Actually, these two classes are families of manifolds whose properties are functions of $h$. The classes regarding $F$ corresponding to $\Bia$(VI$_h$), $h<0$ and $\Bia$(VII$_h$), $h>0$, according to Theorem~4.1, can not be restricted for special values of $h$.

In this section an object of special interest are the curvature properties of these manifolds in relation with $h$.

Having in mind Proposition~3.1, 
it is reasonable to investigate all three subtypes of the Bianchi classes $\Bia(\mathrm{VI}_h)$, $h\leq0$ and $\Bia(\mathrm{VII}_h)$, $h\geq0$.

\subsection{$\Bia(\mathrm{VI}_{\lowercase{h}}),\; \lowercase{h}\leq0$.} 

Let us consider the subtype (1) of this Bianchi class as it is given in Proposition~3.1: 
\[
[e_1,e_2]=o,\quad [e_2,e_3]=e_1-he_2,\quad [e_3,e_1]=he_1-e_2.
\]
The nonzero components of $\n$ for $\Bia(\mathrm{VI}_h$) are calculated:
\begin{equation}\label{VI-nabli}
\begin{array}{lll}
\n_{e_1}e_1=he_3, \quad &
\n_{e_1}e_3=-he_1, \quad &
\n_{e_2}e_2=-he_3, \\
\n_{e_2}e_3=-he_2,  \quad &
\n_{e_3}e_1=-e_2,  \quad &
\n_{e_3}e_2=-e_1.
\end{array}
\end{equation}

Using \eqref{mi:snf}, \eqref{strL}, \eqref{gij} and \eqref{VI-nabli}, we obtain the square norm of $\nabla \f$ as follows
\begin{equation}\label{VI1-norm}
\norm{\n\f}=4(2-h^2).
\end{equation}

Further, there are computed the basic
components $R_{ijkl}=R(e_i,e_j,e_k,e_l)$ of the curvature tensor
$R$, $\rho_{jk}=\rho(e_j,e_k)$ of the
Ricci tensor $\rho$, $\rho^*_{jk}=\rho^*(e_j,e_k)$ of the
associated Ricci tensor $\rho^*$, the values
of the scalar curvatures $\tau$ and $\tau^*$ and of the sectional curvatures $k_{ij}=k(e_i,e_j)$ as
follows:
\begin{equation}\label{VI1-R,rho,tau,k}
\begin{array}{c}
R_{1212}=-R_{1313}=R_{2323}=-h^2;
\\
\rho_{11}=-\rho_{22}=\rho_{33}=-2h^2,\qquad
\rho^*_{12}=\rho^*_{21}=-h^2;
\\
\tau=-6h^2, \qquad \tau^*=0;
\\
k_{12}=k_{13}=k_{23}=-h^2.
\end{array}
\end{equation}

Using the latter equalities we can conclude the following
\vspace*{12pt}

{\bf Proposition~5.1}. 
In the case $\Bia(\mathrm{VI}_h)$, subtype (1), the
following statements are valid:
\begin{enumerate}
\item $\Lf$ is flat if and only if $h=0$;
\item $\Lf$ is an isotropic-cosymplectic B-metric manifold if and only if $h=-\sqrt{2}$;
\item The scalar curvature and the sectional curvatures are constant and non-positive;
\item $\Lf$ is $*$-scalar flat, \ie $\tau^*=0$;
\item $\Lf$ is an Einstein manifold.
\end{enumerate}
\vspace*{12pt}

By similar way we obtain the corresponding results to \eqref{VI1-norm} and \eqref{VI1-R,rho,tau,k} for the rest cases. For the subtype (2) we get:
\[
\begin{array}{c}
\norm{\n\f}=2(1-5h^2);
\\
R_{1212}=-R_{1313}=R_{2323}=-h^2;
\\
\rho_{11}=-\rho_{22}=\rho_{33}=-2h^2,\qquad
\rho^*_{12}=\rho^*_{21}=-h^2;
\\
\tau=-6h^2, \qquad \tau^*=0;
\\
k_{12}=k_{13}=k_{23}=-h^2.
\end{array}
\]

These results imply the following
\vspace*{12pt}

{\bf Proposition~5.2}.
In the case $\Bia$(VI$_h$), subtype (2), the
following statements are valid:
\begin{enumerate}
\item $\Lf$ is flat if and only if $h=0$;
\item $\Lf$ is an isotropic-cosymplectic B-metric manifold if and only if $h=-\frac{\sqrt{5}}{5}$;
\item The scalar curvature and the sectional curvatures are constant and non-positive;
\item $\Lf$ is $*$-scalar flat;
\item $\Lf$ is an Einstein manifold.
\end{enumerate}
\vspace*{12pt}

In the case of subtype (3) we obtain:
\[
\begin{array}{c}
\norm{\n\f}=10(h^2+1);
\\
R_{1212}=R_{2323}=h^2+1,\qquad
R_{1313}=1-h^2,\qquad
R_{1223}=2h;
\\
\rho_{11}=\rho_{33}=2h^2,\qquad
\rho_{13}=\rho_{31}=-2h,\qquad
\rho_{22}=-2(h^2+1);
\\
\rho^*_{12}=\rho^*_{21}=h^2+1,\qquad
\rho^*_{23}=\rho^*_{32}=-2h;\qquad
\\
\tau=2(3h^2+1), \qquad \tau^*=0;
\\
k_{12}=k_{23}=h^2+1,\qquad
k_{13}=h^2-1.
\end{array}
\]

The latter equalities imply the following
\vspace*{12pt}

{\bf Proposition~5.3}.
In the case $\Bia$(VI$_h$), subtype (3), the
following statements are valid:
\begin{enumerate}
\item The square norm of $\n\f$ and the scalar curvature are positive;
\item $\Lf$ is $*$-scalar flat;
\item The sectional curvatures of the $\f$-holomorphic sections are constant and positive.
\end{enumerate}

\subsection{$\Bia(\mathrm{VII}_{\lowercase{h}}),\; \lowercase{h}\geq0$.}

In this subsection we focus our investigations on the three subtypes of $\Bia(\mathrm{VII}_h)$. Firstly, let us consider the subtype (1).
By similar way as the previous subsection, we obtain:
\[
\begin{array}{c}
\norm{\n\f}=4(1-h^2);
\\
R_{1212}=-(h^2+1),\qquad
R_{1313}=-R_{2323}=h^2-1,\qquad
R_{1323}=-2h;
\\
\rho_{11}=-\rho_{22}=-2h^2,\qquad
\rho_{12}=\rho_{21}=2h,\qquad
\rho_{33}=2(1-h^2);
\\
\rho^*_{12}=\rho^*_{21}=-(h^2+1),\qquad
\rho^*_{33}=4h;\qquad
\\
\tau=2(1-3h^2), \qquad \tau^*=4h;
\\
k_{12}=-(h^2+1),\qquad
k_{13}=k_{23}=1-h^2.
\end{array}
\]

These results imply the following
\vspace*{12pt}

{\bf Proposition~5.4}.
In the case $\Bia(\mathrm{VII}_h)$, subtype (1), the
following statements are valid:
\begin{enumerate}
\item $\Lf$ is an isotropic-cosymplectic B-metric manifold if and only if $h=1$;
\item $\Lf$ is scalar flat if and only if $h=\frac{\sqrt{3}}{3}$;
\item $\Lf$ is $*$-scalar flat if and only if $h=0$;
\item The sectional curvatures of the $\f$-holomorphic sections are constant and negative;
\item The sectional curvatures of the $\xi$-sections are constant;
\item $\Lf$ is an $\eta$-complex-Einstein manifold.
\end{enumerate}

Analogously, we get the corresponding results for subtype (2):
\[
\begin{array}{c}
\norm{\n\f}=-10(h^2-1);
\\
R_{1212}=-R_{1313}=-(h^2-1),\qquad
R_{2323}=-(h^2+1),\qquad
R_{1213}=2h;
\\
\rho_{11}=-2(h^2-1),\qquad
\rho_{22}=-\rho_{33}=2h^2,\qquad
\rho_{23}=\rho_{32}=-2h;
\\
\rho^*_{12}=\rho^*_{21}=-(h^2-1),\qquad
\rho^*_{13}=\rho^*_{31}=2h;\qquad
\\
\tau=-2(3h^2-1), \qquad \tau^*=0;
\\
k_{12}=k_{13}=-(h^2-1),\qquad
k_{23}=-(h^2+1).
\end{array}
\]

The latter equalities implies the following
\vspace*{12pt}

{\bf Proposition~5.5}.
In the case $\Bia$(VII$_h$), subtype (2), the
following statements are valid:
\begin{enumerate}
\item $\Lf$ is an isotropic-cosymplectic B-metric manifold if and only if $h=1$;
\item $\Lf$ is scalar flat if and only if $h=\frac{\sqrt{3}}{3}$;
\item $\Lf$ is $*$-scalar flat;
\item $\Lf$ is horizontal flat, \ie  $R|_H=0$ for $H=\ker(\eta)$, if and only if $h=1$;
\item $\rho^*$ and $\tg$ are proportional on $H$ as $\rho^*|_H=(h^2-1)\tg|_H$;
\item $\Lf$ is horizontal $*$-Ricci flat, \ie  $\rho^*|_H=0$, if and only if $h=1$.
\end{enumerate}
\vspace*{12pt}

For the case of the subtype (3) we have:
\[
\begin{array}{c}
\norm{\n\f}=2(5h^2+1);
\\
R_{1212}=-R_{1313}=R_{2323}=h^2;
\\
\rho_{11}=-\rho_{22}=\rho_{33}=2h^2;
\\
\rho^*_{12}=\rho^*_{21}=h^2;
\\
\tau=6h^2, \qquad \tau^*=0;
\\
k_{12}=k_{13}=k_{23}=h^2.
\end{array}
\]

We can conclude the following
\vspace*{12pt}

{\bf Proposition~5.6}.
In the case $\Bia$(VII$_h$), subtype (3), the
following statements are valid:
\begin{enumerate}
\item $\Lf$ is flat if and only if $h=0$;
\item The square norm of $\n\f$ is positive;
\item $\Lf$ is $*$-scalar flat;
\item The scalar curvature and the sectional curvatures are constant and non-negative;
\item $\Lf$ is an Einstein manifold.
\end{enumerate}

\vspace{14pt}

\centerline{REFERENCES}

\parindent 0pt

\small
\begin{enumerate} \frenchspacing

\item{ 
Bianchi, L. 
Sugli spazi a tre dimensioni che ammettono un gruppo
continuo di movimenti. 
{\it Memorie di Matematica e di Fisica della Societa
Italiana delle Scienze, Serie Terza}, {\bf 11}, 1898, 267--352.}

\item{ \vskip-2pt 
Bianchi, L.
 On the three-dimensional spaces which admit a continuous group of motions. 
 {\it Gen. Rel. Grav.}, {\bf 33}, 2001, 2171--2253.}

\item{ \vskip-2pt 
Blair, D. E.
Riemannian Geometry of Contact and Symplectic
Manifolds. 
Pro\-gress in Mathematics, {\bf 203}, Birkh\"auser,
Boston, 2002.}

\item{ \vskip-2pt 
Fagundes, H.
Closed spaces in cosmology. 
{\it Gen. Rel. Grav.},
{\bf 24}, 1992, 199--217.}

\item{ \vskip-2pt 
Ganchev, G., Borisov, A.
Note on the almost complex manifolds with a Norden metric. 
{\it C. R. Acad. Bulg. Sci.}, {\bf 39}, 1986, 31--34.}

\item{ \vskip-2pt 
Ganchev, G., Mihova, V., Gribachev, K.
Almost contact manifolds with B-metric. 
{\it Math. Balkanica (N.S.)},  {\bf 7}, (3-4), 1993, 261--276.}

\item{ \vskip-2pt 
Gribachev, K., Mekerov, D., Djelepov, G.
Generalized B-manifolds. 
{\it C. R. Acad. Bulg. Sci.}, {\bf 38}, 1985, 299--302.}

\item{ \vskip-2pt 
Manev, H.
On the structure tensors of almost contact
B-metric manifolds. arXiv:\allowbreak{}1405.3088.}

\item{ \vskip-2pt 
Manev, H., Mekerov, D.
Lie groups as 3-dimensional almost contact B-metric manifolds. 
{\it J. Geom.}, Online First (17 Sept 2014), DOI 10.1007/s00022-014-0244-0.}

\item{ \vskip-2pt 
Manev, M.
Natural connection with totally skew-symmetric
tor\-sion on almost con\-tact manifolds with B-metric. 
{\it Int. J. Ge\-om. Methods Mod. Phys.}, {\bf 9}, (5), 2012, 1250044 (20 pa\-ges).}

\item{ \vskip-2pt 
Manev, M., Nakova, G.
Curvature properties of some
three-dimensional almost con\-tact B-metric manifolds. 
{\it Plovdiv Univ. Sci. Works -- Math.}, {\bf{34}}, (3), 2004, 51--60.}

\item{ \vskip-2pt 
Nakova, G., Manev, M.
Curvature properties on some
three-dimensional almost con\-tact manifolds with B-metric, II.
In {\em Proc. 5th Internat.
Conf. Geometry, Integrability \& Qu\-antization\/} V, eds. I. M. Mladenov and A. C. Hirshfeld, 2004, 169--177.}

\item{ \vskip-2pt 
Thurston, W.
Three-dimensional geometry and topology. Vol. 1. Princeton University Press, 1997.}

\end{enumerate}
  \vspace{12pt}
\baselineskip10pt

\vskip10pt
\begin{flushright}
\end{flushright}
\vskip20pt
 \footnotesize
\begin{flushleft}
Department of Pharmaceutical Sciences\\
Faculty of Pharmacy\\
Medical University of Plovdiv\\
  15-A  blvd. Vasil Aprilov,   BG-4002 Plovdiv\\
BULGARIA \\
and\\
Department of Algebra and Geometry\\
Faculty of Mathematics and Informatics\\
Paisii Hilendarski University of Plovdiv\\
236  blvd. Bulgaria, BG-4027 Plovdiv\\
BULGARIA \\
e-mail: hmanev@uni-plovdiv.bg
\end{flushleft}

\end{document}